\documentclass[12pt]{article}

\usepackage{amsfonts}
\usepackage{amssymb}
\usepackage{amsmath,bm}
\usepackage{fancyhdr}    
\usepackage{epic}
\usepackage{eepic}
\usepackage{epsf}
\usepackage{graphicx}
\usepackage{rotating}

\DeclareMathAlphabet{\pazocal}{OMS}{zplm}{m}{n}
\def\dia{~~$\diamondsuit$}

\def\bea{\begin{eqnarray}}
\def\eea{\end{eqnarray}}

\def\eqnn#1{\bea\label{#1}}

 \renewcommand{\leq}{\leqslant}
\renewcommand{\geq}{\geqslant}

\def\eps{\epsilon}

\def\til{{\tilde I^\L}}
\def\bd{{\bar \D}}

\def\nd{\end{document}}

\def\qrt{{\textstyle{1\over4}}}
\def\ha{{\textstyle{1\over2}}}

\def\({\left(}\def\){\right)}
\def\hH{{\hat H}} 
\def\k{\kappa}
\def\rg{\rangle}

\def\r{\rho}
\def\half{{\textstyle{\frac{1}{2}}}}
\def\trh{{\textstyle{\frac{3}{2}}}}
\def\frh{{\textstyle{\frac{5}{2}}}}

\def\riga{-\kern-5pt - \kern-5pt -}
\def\white#1{{\bigcirc\atop {#1}}}
\def\black#1{{\bullet \atop {#1}}}

\def\third{{\textstyle{1\over3}}}

\def\o{{\bar 0}} \def\I{{\bar 1}}
\def\a{\alpha}
\def\b{\beta}
\def\d{\delta}

\def\D{\Delta} \def\L{\Lambda}

\def\bbz{\mathbb{Z}}
\def\bbc{\mathbb{C}}
\def\bac{\bbc} 
\def\bbr{\mathbb{R}}
\def\bbn{\mathbb{N}}

\def\be{\begin{equation*}}
\newcommand{\ees}{\end{equation*}}

\def\eqn{\begin{equation}\label}
\newcommand{\ee}{\end{equation}}

\def\beq{\begin{eqnarray*}}
\def\eeq{\end{eqnarray*}}

\def\bea {\begin {eqnarray}}
\def\eqnn#1{\bea\label{#1}}
\def\eea{\end{eqnarray}}

\def\nn{\nonumber}

\newcommand{\eqna}[1]{\begin{subequations} \label{#1}
\begin{eqnarray}}
\def\eena{\end{eqnarray}
\end{subequations}}

\def\vr{\vert} \def\bbr{I\!\!R}

\def\nt{\noindent}

\def\cg{{\pazocal G}} \def\ch{{\pazocal H}} 
 \def\ck{{\pazocal K}} 
  
\def\cp{{\pazocal P}}

\newcommand{\ervv}[1]{{Y_{#1}}} 
\newcommand{\orv}[1]{{Y_{#1}}} 
\newcommand{\lwv}{v_0} 

\input epsf.tex
\newcount\figno
\def\fig#1#2#3{
\par\begingroup\parindent=0pt\leftskip=1cm\rightskip=1cm\parindent=0pt
\baselineskip=11pt \global\advance\figno by 1 
\epsfxsize=#3 \centerline{\epsfbox{#2}} \vskip 12pt
#1\par
\endgroup\par}
\def\figlabel#1{\xdef#1{\the\figno}}
\def\encadremath#1{\vbox{\hrule\hbox{\vrule\kern8pt\vbox{\kern8pt
\hbox{$\displaystyle #1$}\kern8pt} \kern8pt\vrule}\hrule}}

\begin{document}
\pagestyle{empty}

\begin{center}

\textsf{\Large\bf Positive Energy Unitary Irreducible\\[3pt]
Representations of the Superalgebra \bm{$osp(1|8,\bbr)$}\footnote{To appear in
Publications de l'Institut Math\'ematique,
Nouvelle s\'erie, 
(2017)}}


\vspace{7mm}

{\large  V.K. Dobrev$^{a,}$\footnote{dobrev@inrne.bas.bg, Supported partially by COST Actions MP1210 and MP1405, and by  Bulgarian NSF Grant DFNI T02/6.} and I.
Salom$^{b,}$\footnote{isalom@phy.bg.ac.rs, Supported partially by COST Action MP1405 and Serbian Ministry of Science and
Technological Development, grant OI 171031}}

\vspace{3mm}

 \emph{$^a$Institute of Nuclear Research and Nuclear Energy,\\
 Bulgarian Academy of Sciences, \\
72 Tsarigradsko Chaussee, 1784 Sofia, Bulgaria}

\vspace{2mm}

 \emph{$^b$Institute of Physics, University of
Belgrade,\\ Pregrevica 118, 11080 Zemun, Belgrade, Serbia}

\end{center}

\vspace{.8 cm}



\begin{abstract}
We continue the study of positive energy (lowest weight) unitary
irreducible representations of the superalgebras ~$osp(1|2n,\bbr)$.
We present the full list of these UIRs. We give the Proof of the case
~$osp(1|8,\bbr)$.
\end{abstract}

\section{Introduction}

Recently, superconformal field theories in various dimensions are
attracting more interest, in particular, due to their duality to AdS
supergravities. This makes the classification of the UIRs of these
superalgebras very important. Until recently only those for ~$D\leq
6$~ were studied since in these cases the relevant superconformal
algebras satisfy \cite{Nahm} the Haag-Lopuszanski-Sohnius theorem
\cite{HLS}. Thus, such classification was known only for the ~$D=4$~
superconformal algebras ~$su(2,2/N)$ \cite{FF} (for $N=1$),
\cite{DPm,DPu,DPf,DPp} (for arbitrary $N$). More recently, the
classification for ~$D=3$ (for even $N$), $D=5$, and $D=6$ (for
$N=1,2$)~ was given in \cite{Min} (some results are conjectural),
and then the $D=6$ case (for arbitrary $N$) was finalized in
\cite{Dosix}.

On the other hand the applications in string theory require the
knowledge of the UIRs of the conformal superalgebras for ~$D>6$.
Most prominent role play the superalgebras $osp(1\vr\,2n)$.
Initially, the superalgebra $osp(1\vr\,32)$ was put forward for
$D=10$ \cite{Tow}. Later it was realized that $osp(1\vr\,2n)$ would
fit any dimension, though they are minimal only for $D=3,9,10,11$
(for $n=2,16,16,32$, resp.) \cite{DFLV}. In all cases we need to
find first the UIRs of ~$osp(1\vr\, 2n,\bbr)$~ which study was
started in \cite{DZ} and \cite{DMZZ}. Later, in \cite{DS} we
finalized the UIR classification of \cite{DZ} as Dobrev-Zhang-Salom
(DZS) Theorem. In \cite{DS} we proved the DZS Theorem for
$osp(1\vr\,6)$.

In the present paper we prove the DZS Theorem for $osp(1\vr\,8)$.
For the lack of space we refer for extensive literature on the
subject in \cite{DZ,DS}.

\section{Preliminaries on representations}

 Our basic references for Lie superalgebras are \cite{Kab,Kc},
although in this exposition we follow \cite{DZ}.

The even subalgebra of ~$\cg = osp(1\vr\, 2n,\bbr)$~ is the algebra
~$sp(2n,\bbr)$ with maximal compact subalgebra ~$\ck = u(n) \cong
su(n) \oplus u(1)$.

We label the relevant representations of ~$\cg$~
by the signature:
\eqn{sgn}\chi ~=~ [\, d\,;\,a_1\,,...,a_{n-1}\,] \ee
where ~$d$~ is the conformal weight, and ~$a_1,...,a_{n-1}$~
are non-negative integers which are Dynkin labels of the
finite-dimensional UIRs of the subalgebra $su(n)$ (the simple
part of $\ck$).

In \cite{DZ} were classified (with some omissions to be spelled out
below) the positive energy (lowest weight) UIRs of ~$\cg$~ following
the methods used for the $D=4,6$ conformal superalgebras, cf.
\cite{DPm,DPu,DPf,DPp,Dosix}, resp. The main tool was an adaptation
of the Shapovalov form \cite{Sha} on the Verma modules ~$V^\chi$~
over the complexification ~$\cg^\bac ~=~ osp(1\vr\, 2n)$~ of ~$\cg$.

We recall some facts about ~$\cg^\bac ~=~ osp(1\vr\, 2n)$
(denoted $B(0,n)$ in \cite{Kab}) as used in \cite{DZ}.
The root systems are given in terms of
~$\d_1\,\dots,\d_{n}\,$, ~$(\d_i,\d_j) ~=~
\d_{ij}\,$, ~$i,j=1,...,n$.
The even and odd roots systems are \cite{Kab}:
\beq
\D_\o &=& \{ \pm\d_i\pm\d_j\ , ~1\leq i< j\leq n\ ,
~~\pm 2\d_i\ , ~1\leq i \leq n \} ~,
\\ \D_\I &=& \{ \pm\d_i\ , ~1\leq i \leq n
\} \eeq
(we remind that the signs ~$\pm$~ are not correlated).
We shall use the following distinguished
simple root system \cite{Kab}:
\be
 \Pi ~=~ \{\, \d_1-\d_2 \, ,\,
 , \dots, \d_{n-1}- \d_{n} \, ,\, \d_n\, \} \ , \ees
or introducing standard notation for the simple roots:
\beq
\Pi ~&=&~ \{\,\a_1\,,...,\,\a_{n} \, \}\ , \\
&&\a_{j} ~=~ \d_{j}-\d_{j+1}\ , \quad
j=1,...,n-1 \ , \quad \a_{n} ~=~ \d_{n}\ . \eeq
The root ~$\a_n = \d_n$~ is odd, the other simple roots
are even.
The Dynkin diagram is:{\Large
\be
\vbox{\offinterlineskip\baselineskip=10pt
\halign{\strut#
\hfil
& #\hfil
\cr &\cr
&\ $\white{1}\
{\riga \cdots \riga\atop \phantom{a}}
\white{{n-1}} {=\kern-2pt\Longrightarrow \atop \phantom{a}}
\black{n}$
\cr }} \ees}
The black dot is used to signify that the simple odd root is not
nilpotent. In fact, the
superalgebras ~$B(0,n) = osp(1\vr\, 2n)$~ have no nilpotent generators
unlike all other types of basic classical Lie superalgebras \cite{Kab}.

The corresponding to ~$\Pi$~ positive root system is:
\eqn{psdm}\D_\o^+ \,=\, \{ \d_i\pm\d_j\,, \ 1\leq i< j\leq n,
 \,\,2\d_i\,, \ 1\leq i \leq n \} ,
\qquad \D_\I^+ \,=\, \{ \d_i \,, \ 1\leq i \leq n \} \ee
We record how the elementary functionals are expressed through
the simple roots:
\be \d_k ~=~ \a_k + \cdots + \a_n \ .\ees

From the point of view of representation theory more relevant is the restricted root system, such that:
\beq && \bd^+ ~=~ \bd^+_\o \cup \D_\I^+ \ , \\
&& \bd^+_\o ~\equiv ~ \{ \a\in\D^+_\o ~\vert ~\half\a\notin\D^+_\I
\} ~=~ \{ \d_i\pm\d_j\ , ~1\leq i< j\leq n \} \eeq

The superalgebra ~$\cg ~=~ osp(1\vr\, 2n,\bbr)$~ is a split real form
of $osp(1\vr\, 2n)$ and has the same root system.

The above simple root system is also the simple root system of the complex
simple Lie algebra ~$B_n$ (dropping the distinction between even and odd roots) with
Dynkin diagram:{\Large
\be
\vbox{\offinterlineskip\baselineskip=10pt
\halign{\strut#
\hfil
& #\hfil
\cr &\cr
&\ $\white{{1}}
{\riga \cdots \riga\atop \phantom{a}}
\white{{n-1}} {=\kern-2pt\Longrightarrow\atop \phantom{a}}
\white{n}$
\cr }} \ees}
Naturally, for the $B_n$ positive root system we drop the roots $2\d_i\,$~
\be\D^+_{\rm B_n} ~=~ \{ \d_i\pm\d_j\ , ~1\leq i< j\leq n\ ,
 ~~\d_i \ , ~1\leq i \leq n \} ~\cong ~\bd^+ \ees
This shall be used essentially below.

Besides \eqref{sgn} we shall use the Dynkin-related labelling:
\be
(\L , \a_k^\vee ) ~= -\, a_k \ , ~~~ 1\leq k \leq n\ , \ees
where ~$\a_k^\vee \equiv 2\a_k/(\a_k,\a_k)$, and
the minus signs
are related to the fact that we work with
lowest weight Verma modules (instead of the highest weight modules
used in \cite{Kc}{}) and to Verma module reducibility
w.r.t. the roots ~$\a_k$ (this is explained in detail in \cite{DPf,DZ}).

Obviously, ~$a_n$~ must be related to the conformal weight ~$d$~ which is
 a matter of
normalization so as to correspond to some known cases.
Thus, our choice is:
\be a_n ~=~ - 2 d - a_1 - \cdots - a_{n-1} \ . \ees

 The actual Dynkin labelling is given by:
\be m_k ~=~ (\r - \L, \a_k^\vee )\ees
where ~$\r\in\ch^*$~ is given by the
difference of the half-sums $\r_\o\,,\r_\I$
of the even, odd, resp., positive roots (cf. (\ref{psdm}):
\beq \r ~&\doteq&~ \r_\o - \r_\I
~=~ (n-\half)\d_1 + (n-\trh) \d_2 + \cdots + \trh \d_{n-1} +
\half \d_n \ , \\
&&\r_\o\ =\ n \d_1 + (n-1) \d_2 + \cdots +
2\d_{n-1} + \d_n \ , \nn\\
&&\r_\I \ =\ \half (\d_1 + \cdots + \d_n)\ . \nn\eeq

Naturally, the value of $\r$ on the simple roots is 1:
~$(\r,\a^\vee_i)=1$, $i=1,...,n$.

Unlike $a_k\in\bbz_+$ for $k<n$ the value of $a_n$ is arbitrary.
In the cases when $a_n$ is also a non-negative integer, and then ~$m_k\in\bbn$ ($\forall k$) the
corresponding irreps are the finite-dimensional irreps of $\cg$ (and of $B_n$).

Having in hand the values of ~$\L$~ on the basis we can
recover them for any element of ~$\ch^*$.
We shall need only ~$(\L , \b^\vee)$~ for all positive roots $\b$
 as given in \cite{DZ}:
\eqnn{lorr}
(\L , (\d_i-\d_j)^\vee ) ~&=&~ (\L , \d_i-\d_j ) ~=~
- a_i - \cdots - a_{j-1} \nn\\
(\L , (\d_i+\d_j)^\vee ) ~&=&~ (\L , \d_i+\d_j ) ~=~
2d\ +\ a_1 + \cdots + a_{i-1} - a_j - \cdots - a_{n-1} \qquad
\nn\\
(\L , \d_i^\vee ) ~&=&~ (\L , 2\d_i ) ~=~
2d \ +\ a_1 + \cdots + a_{i-1} - a_i - \cdots - a_{n-1}
\\
(\L , (2\d_i)^\vee ) ~&=&~ (\L , \d_i ) ~=~
d \ +\ \half ( a_1 + \cdots + a_{i-1} - a_i - \cdots - a_{n-1} ) \nn
\eea

To introduce Verma modules we use the standard
triangular decomposition:
\be \cg^\bac ~=~ \cg^+ \oplus \ch \oplus \cg^- \ees
where $\cg^+$, $\cg^-$, resp., are the subalgebras corresponding
to the positive, negative, roots, resp., and $\ch$ denotes the
Cartan subalgebra.

We consider lowest weight Verma modules,
so that ~$V^\L ~ \cong U(\cg^+) \otimes v_0\,$,\\
 where ~$U(\cg^+)$~ is the universal enveloping algebra of $\cg^+$,
and ~$v_0$~ is a lowest weight vector $v_0$ such that:
\beq
 Z \ v_0\ &=&\ 0 \ , \quad Z\in \cg^- \nn\\
 H \ v_0 \ &=&\ \L(H)\ v_0 \ , \quad H\in \ch \ .\eeq
Further, for simplicity we omit the sign
~$\otimes \,$, i.e., we write $p\,v_0\in V^\L$ with $p\in U(\cg^+)$.

Adapting the criterion of \cite{Kc} (which generalizes the
BGG-criterion \cite{BGG} to the super case) to lowest weight
modules, one finds that a Verma module ~$V^\L$~ is reducible w.r.t.
the positive root ~$\b$~ iff the following holds \cite{DZ}:
\eqn{odr} (\r - \L, \b^\vee) = m_\b\ , \qquad \b \in \Delta^+ \ ,
\quad m_\b\in\bbn \ . \ee

If a condition from (\ref{odr}) is fulfilled then ~$V^\L$~ contains
a submodule which is a Verma module ~$V^{\L'}$~ with shifted weight
given by the pair ~$m,\b$~: ~$\L' ~=~ \L + m\b$. The embedding of
~$V^{\L'}$~ in ~$V^\L$~ is provided by mapping the lowest weight
vector ~$v'_0$~ of ~$V^{\L'}$~ to the ~{\bf singular vector}
~$v_s^{m,\b}$~ in ~$V^\L$~ which is completely determined by the
conditions:
\beq
 X \ v_s^{m,\b}\ &=&\ 0 \ , \quad X\in \cg^- \ , \nn\\
 H \ v_s^{m,\b} \ &=&\ \L'(H)\ v_0 \ , \quad H\in \ch \ ,
~~~\L' ~=~ \L + m\b\ .\eeq
Explicitly, ~$v_s^{m,\b}$~ is given by a
polynomial in the positive root generators \cite{Dob,DPf}:
\be v_s^{m,\b} ~=~ P^{m,\b} \,v_0 \ , \quad P^{m,\b}\in
U(\cg^+)\ . \ees
Thus, the submodule ~$I^\b$~ of ~$V^\L$~ which is
isomorphic to ~$V^{\L'}$~ is given by\\ ~$U(\cg^+)\, P^{m,\b}
\,v_0\,$.

Note that the Casimirs of $\cg^\bac$ take the same values on
~$V^\L$~ and ~$V^{\L'}$.

Certainly, \eqref{odr} may be fulfilled for several positive roots
(even for all of them). Let ~$\D_\L$~ denote the set of all positive
roots for which \eqref{odr} is fulfilled, and let us denote:
~$\tilde I^\L ~\equiv~ \cup_{\b\in\D_\L} I^{\b}\,$. Clearly,
~$\tilde I^\L$~ is a proper submodule of ~$V^\L$. Let us also denote
~$F^\L ~\equiv V^\L / \tilde I^\L$.

Further we shall use also the following notion. The singular vector
~$v_1$~ is called ~{\bf descendant}~ of the singular vector
~$v_2\notin \bbc v_1$~ if there exists a homogeneous polynomial
~$P_{12}$~ in ~$U(\cg^+)$~ so that ~$v_1 ~=~ P_{12} ~v_2\,$. Clearly,
in this case we have: ~$I^1 ~\subset ~I^2\,$, where ~$I^k$~ is the
submodule generated by $v_k\,$.

The Verma module ~$V^\L$~ contains a unique proper maximal submodule
~$I^\L$ ($\supseteq \til$) \cite{Kc,BGG}. Among the lowest weight
modules with lowest weight ~$\L$~ there is a unique irreducible one,
denoted by ~$L_\L$, i.e., ~$L_\L ~=~ V^\L/I^\L$. (If ~$V^{\L}$~ is
irreducible then ~$L_\L = V^\L$.)

It may happen that the maximal submodule ~$I^\L$~ coincides with
the submodule ~$\tilde I^\L$~ generated by all singular vectors.
This is, e.g., the case for all Verma modules if ~rank~$\cg ~\leq
2$, or when \eqref{odr} is fulfilled for all simple roots (and, as a
consequence for all positive roots). Here we are interested in the
cases when ~$\tilde I^\L$~ is a proper submodule of ~$I^\L$. We need
the following notion.

{\bf Definition:} \cite{BGG,Docond,Do-KL} ~{\it Let ~$V^\L$~ be a
reducible Verma module. A vector ~$v_{{\rm ssv}} ~\in ~V^\L$~ is
called a ~{\bf subsingular vector}~ if ~$v_{{\rm su}} ~\notin ~
\tilde I^\L$~ and the following holds:}
\be X~v_{{\rm su}}
~~\in ~~ \tilde I^\L \,, \quad \forall X\in\cg^- \ees

\medskip

Going from the above more general definitions to ~$\cg$~ we recall
that in \cite{DZ} it was established that from \eqref{odr} follows
that the Verma module ~$V^{\L(\chi)}$~ is reducible if one of the
following relations holds (following the order of (\ref{lorr}):
\eqna{redd}
&\bbn \ni m^-_{ij} = j-i + a_i + \cdots + a_{j-1} \\
&\bbn \ni m^+_{ij} = 2n -i-j +1 + a_j + \cdots + a_{n-1}
- a_1 - \cdots - a_{i-1} -2d \\
&\bbn \ni m_i = 2n-2i+1 + a_i + \cdots + a_{n-1} - a_1 +
\cdots - a_{i-1} - 2d \\
&\bbn \ni m_{ii} = n-i+\half(1 + a_i + \cdots + a_{n-1} - a_1 +
\cdots - a_{i-1}) - d \ . \eena
Further we shall use the fact from \cite{DZ} that we may eliminate the reducibilities and
embeddings related to the roots ~$2\d_i\,$. Indeed,
since $m_i ~=~ 2m_{ii}\,$, whenever (\ref{redd}d) is
fulfilled also (\ref{redd}c) is fulfilled.

For further use we introduce notation for the root vector ~$X^+_j\,\in\cg^+$, ~$j=1,\ldots,n$,
~corresponding to the simple root ~$\a_j\,$. Naturally, ~$X^-_j\,\in\cg^-$~ corresponds to ~$-\a_j\,$.

Further, we notice that all reducibility conditions in
(\ref{redd}{a}) are fulfilled. In particular, for the simple roots
from those condition (\ref{redd}{a}) is fulfilled with
~$\b\to\a_i=\d_i-\d_{i+1}\,$, $i=1,...,n-1$ and ~$m^-_i ~\equiv~
m^-_{i,i+1} ~=~ 1 + a_i\,$. The corresponding submodules ~$I^{\a_i}
~=~ U(\cg^+)\,v^i_s\,$, where ~$\L_i ~=~ \L + m^-_i \a_i$~ and
~$v^i_s ~=~ (X^+_i)^{1+a_i}\, v_0\,$. These submodules generate an
invariant submodule which we denote by ~$I^\L_c\, \subset \tilde
I^\L$. Since these submodules are nontrivial for all our signatures
in the question of unitarity instead of ~$V^\L$~ we shall consider
also the factor-modules:
\be F_c^\L ~=~ V^\L\, /\, I^\L_c ~\supset ~F^\L\ . \ees
We shall denote the lowest weight vector of
~$F_c^\L$~ by ~$\vr \L_c\rg$~ and the singular vectors above become
null conditions in~ $F_c^\L$~:
\be (X^+_i)^{1+a_i}\, \vr \L_c\rg ~=~ 0 \ , \quad i=1,...,n-1. \ees

If the Verma module ~$V^\L$~ is not reducible w.r.t. the other
roots, i.e., (\ref{redd}{b,c,d}) are not fulfilled, then
~$F_c^\L=F^\L$~ is irreducible and is isomorphic to the irrep
~$L_\L$~ with this weight.

In fact, for the factor-modules reducibility is controlled by the value of ~$d$,
or in more detail:

The maximal ~$d$~ coming from the different possibilities in
(\ref{redd}{b}) are obtained for $m^+_{ij}=1$ and they are:
\be
d_{ij} ~\equiv~ n + \half(a_j + \cdots + a_{n-1}
- a_1 - \cdots - a_{i-1}-i-j) \ ,
\ees
the corresponding root being ~$\d_i+\d_j\,$.

The maximal ~$d$~ coming from the different possibilities in
(\ref{redd}{c,d}), resp., are obtained for $m_{i}=1$, $m_{ii}=1$,
resp., and they are:
\beq
&&d_{i} ~\equiv~ n -i + \half(a_i + \cdots + a_{n-1}
- a_1 - \cdots - a_{i-1}) \ ,
\\
&&d_{ii} ~=~ d_i - \half \ ,
\nn\eeq
the corresponding roots being ~$\d_i\,,2\d_j\,$, resp.

There are some orderings between these maximal reduction points \cite{DZ}:
\eqnn{compz}
d_1 ~&>&~ d_2 ~>~ \cdots ~>~ d_n \ , \\
d_{i,i+1} ~&>&~ d_{i,i+2} ~>~ \cdots ~>~ d_{in} \ , \nn\\
d_{1,j} ~&>&~ d_{2,j} ~>~ \cdots ~>~ d_{j-1,j} \ , \nn\\
d_i ~&>&~ d_{jk} ~>~ d_\ell
\ , \qquad i\leq j <k \leq \ell \ . \nn\eea

Obviously the first reduction point is:
\be d_{1} ~=~ n -1 + \half(a_1 + \cdots + a_{n-1}) \ .
\ees

\section{Unitarity}

The first results on the unitarity were given in \cite{DZ}, and then
improved in \cite{DS}. Thus, the statement below should be
called~ {\it Dobrev-Zhang-Salom Theorem}:\smallskip

\nt {\bf Theorem DZS:} ~~ {\it All positive energy unitary irreducible
representations of the superalgebras ~$osp(1\vr\, 2n,\bbr)$~
characterized by the signature ~$\chi$~ in (\ref{sgn}) are obtained
for real ~$d$~ and are given as follows:
\beq &&d ~\geq~ n -1 +
\half(a_1 + \cdots + a_{n-1}) ~=~ d_{1}
\ , \quad a_1 \neq 0 \ , \nn\\
&&d ~\geq~ n - \trh + \half(a_2 + \cdots + a_{n-1} ) ~=~ d_{12}
\ , \quad a_1 = 0,\ a_2\neq 0\ , \nn\\
&&d ~=~ n - 2 + \half(a_2 + \cdots + a_{n-1} ) ~=~d_{2} > d_{13}
\ , \quad a_1 = 0,\ a_2\neq 0\ , \nn\\
&&d ~\geq~ n - 2 + \half(a_3 + \cdots + a_{n-1} ) ~=~ d_{2} ~=~ d_{13}
\ , \quad a_1 = a_2=0,\ a_3\neq 0\ , \nn\\
&&d ~=~ n - \frh + \half(a_3 + \cdots + a_{n-1} ) ~=~ d_{23} ~>~ d_{14}
\ , \quad a_1 = a_2=0,\ a_3\neq 0\ , \nn\\
&&d ~=~ n - 3 + \half(a_3 + \cdots + a_{n-1} ) ~=~ d_3 ~=~ d_{24} ~>~ d_{15}
\ , \quad a_1 = a_2=0,\ a_3\neq 0\ , \nn\\
&& ... \nn\\
&&d ~\geq~ n - 1 -\k + \half(a_{2\k+1} + \cdots + a_{n-1})
\ , \quad a_1 = ... = a_{2\k}=0,\ a_{2\k+1}\neq 0\ , \nn\\
&& \qquad\qquad \k = \half, 1, ..., \half(n-1) \ , \nn\\
&&d ~=~ n - \trh - \k + \half(a_{2\k+1} + \cdots + a_{n-1})
\ , \quad a_1 = ... = a_{2\k}=0,\ a_{2\k+1}\neq 0\ , \nn\\
&& ... \nn\\
&&d ~=~ n - 1-2\k + \half(a_{2\k+1} + \cdots + a_{n-1})
\ , \quad a_1 = ... = a_{2\k}=0,\ a_{2\k+1}\neq 0\ , \nn\\
&& ... \nn\\
 &&d ~\geq~ \half (n-1)\ , \quad a_1 ~=~ ... ~=~ a_{n-1} ~=~ 0 \nn\\
&&d ~=~ \half (n-2)\ , \quad a_1 ~=~ ... ~=~ a_{n-1} ~=~ 0 \nn\\
&& ... \nn\\
&&d ~=~ \half \ , \quad a_1 ~=~ ... ~=~ a_{n-1} ~=~ 0 \nn\\
&&d ~=~ 0 \ , \quad a_1 ~=~ ... ~=~ a_{n-1} ~=~ 0 \nn
 \eeq
 where the last case is the trivial one-dimensional irrep.} \\
Parts of the {\it Proof} were given in \cite{DZ}, while in
\cite{DS} was given a detailed sketch of the Proof. In \cite{DS} was
given also the Proof for the case $n=3$.

\medskip

In the present paper we give the Proof for $osp(1|8)$.

\bigskip

\setcounter{equation}{0}
\section{The case of osp(1$\vr$8)}

For ~$n=4$~ formula (\ref{compz}) simplifies to:

\fig{}{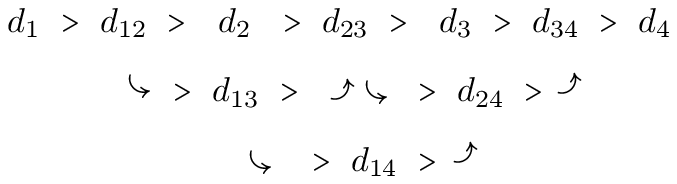}{8cm}


In the case of $osp(1|8)$ the DZS Theorem  reads:\\
\nt {\bf Theorem:}~~
 All positive energy unitary irreducible
representations of the superalgebras ~$osp(1\vr\, 8,\bbr)$~
characterized by the signature ~$\chi$~ in (\ref{sgn}) are obtained
for real ~$d$~ and are given as follows: \beq &&d ~\geq~ 3 +
\half(a_1 + a_{2}+a_3) ~=~d_{1}
\ , \quad a_1\neq 0 \ , \\
&&d ~\geq~ \frh + \half (a_2+a_3) ~=~d_{12}
\ , \quad a_1 = 0,\ a_2\neq 0\ , \nn\\
&&d ~=~ 2 + \half (a_2+a_3) ~=~ d_{2} > d_{13}
\ , \quad a_1 = 0,\ a_2\neq 0\ , \nn\\
&&d ~\geq~ 2 + \half a_3 ~=~ d_2 = d_{13} \ , \quad a_1 = a_2= 0\ , \ a_3\neq 0 \nn\\
&&d ~=~ \trh + \half a_3 ~=~ d_{23} > d_{14} \ , \quad a_1 = a_2=0\ , \ a_3\neq 0 \nn\\
&&d ~=~ 1 + \half a_3 ~=~ d_3 > d_{24} \ , \quad a_1 = a_2= 0\ , \ a_3\neq 0 \nn\\
&&d ~\geq~ \trh ~=~ d_{23} = d_{14} \ , \quad a_1 = a_2= a_3=0 \nn\\
&&d ~=~ 1 ~=~ d_3 = d_{24} \ , \quad a_1 = a_2= a_3= 0 \nn\\
 &&d ~=~ \half ~=~d_{34} \ , \quad a_1 = a_{2}= a_3 = 0 \ , \nn\\
&&d ~=~ 0 ~=~ d_{4} \ , \quad a_1 = a_{2}= a_3 = 0 \ \nn \eeq
where the last case is the trivial one-dimensional irrep.\\
{\bf Proof:}~~
For ~$d>d_1$~ there are no singular vectors and we have unitarity.
At ~$d=d_1$~ there is a singular vector of weight ~$\d_1=\a_1+\a_2+\a_3+\a_4$~ \cite{Dos,DZ}:
\eqnn{svn3d1}
v^1_{\d_1}
\, &=&\, \sum_{k_1=0}^1 \sum_{k_2=0}^1 \sum_{k_{3}=0}^1\, b_{k_1,k_2,k_3}\,
(X^+_1)^{1-k_1} (X^+_2)^{1-k_2} (X^+_{3})^{1-k_{3}}\,
\times \\ &&\times\
X^+_4\, (X^+_{3})^{k_{3}} \, (X^+_{2})^{k_{2}} \,
(X^+_1)^{k_1}\, v_0\, \equiv\, \cp^{1,\d_1}\ v_0 \ ,
\nn\\ \nn\\
 b_{k_1,k_2, k_{3}} \, &=&\, (-1)^{k_1 + k_2 + k_{3}}
\, (a_1+k_1) \, \frac{2+a_1+a_2}{1+a_1+a_2 - k_2} \, \frac{3 +a_1+
a_2 +a_{3}}{3+a_1 + a_2 +a_{3} - k_{3}} \nn\eea
where \ $H^s\ =\ \hH_1 + \hH_2 + \cdots + \hH_s\,$, and a basis in terms of simple root vectors only is used.
This singular vector is non-trivial for ~$a_1\neq 0$~ and must be
eliminated to obtain an UIR. Below ~$d<d_1$~ this vector is not
singular but has negative norm and thus there is no unitarity for
~$a_1\neq 0$. On the other hand for $a_1= 0$ and any $d$ the vector
\eqref{svn3d1} is
descendant of the compact root singular vector ~$X_1^+\,v_0$~ which is already factored out for ~$a_1 = 0$.\\
Thus, below we discuss only the cases with ~$a_1 = 0$~ in which case we have unitarity for
~$d>d_{12} = \frh + \half (a_2+a_3)$.
Then at the next reducibility point ~$d ~=~ d_{12}$~ we have a
 singular vector corresponding to the root ~$\d_1 + \d_2 = \a_1 + 2\a_2 + 2\a_3 +2\a_4$~
which is given by:
\eqnn{vs12}
v^{1}_{\d_1 + \d_2} &=& \frac 1{2 + 2 a_2 + a_3} \Big(
-\frac{1}{2} \left(Y_4 Y_3 X^+_3 (X^+_2 )^2 X^+_1\right)-\frac{1}{4}
\left( Y^2_4 (X^+_3 )^2 (X^+_2 )^2 X^+_1\right)\ + \nn\\  &&+\ \left( Y^2_4 X^+_3
X^+_{23} X^+_2 X^+_1\right) a_2-2 \left(Y_4 Y_2 X^+_{23}
X^+_1\right) a_2 \left(a_2+1\right)\ - \nn\\  &&-\ \left( Y^2_4 X^+_{23} X^+_{23}
X^+_1\right) a_2 \left(a_2+1\right)- \left(Y_4 Y_3 X^+_{23} X^+_2
X^+_1\right) \left(a_3+2\right)
\ - \nn\\  &&-\ 2 \left(Y_3 Y_2 X^+_2 X^+_1\right)
\left(a_2+a_3+1\right) \left(a_2+a_3+2\right)\ - \nn\\  &&-\ \left( Y^2_3  (X^+_2
)^2 X^+_1\right) \left(a_2+a_3+1\right) \left(a_2+a_3+2\right)
\ - \nn\\  &&-\ 4 \left( Y^2_2 X^+_1\right) a_2
\left(a_2+1\right) \left(a_2+a_3+1\right)
\left(a_2+a_3+2\right) \ + \nn\\  &&+\ \left(\ervv {23} X^+_2 X^+_1\right) \left(2
a_2+1\right) \left(a_2+a_3+1\right) \left(a_2+a_3+2\right)
\ + \nn\\  &&+\ \left(Y_4
Y_2 X^+_3 X^+_2 X^+_1\right) \left(2 a_2+a_3+2\right)
\ + \nn\\  &&+\ \qrt
\left(\ervv {34} X^+_3 (X^+_2 )^2 X^+_1\right) \left(2 a_2+2
a_3+3\right)  \ + \nn\\  &&+\ \left(\ervv {24} X^+_{23} X^+_1\right) a_2
\left(a_2+1\right) \left(2 a_2+2 a_3+3\right)
\ + \nn\\  &&+\ \ha
\left(\ervv {34} X^+_{23} X^+_2 X^+_1\right) \left(a_3-2 a_2
\left(a_2+a_3+1\right)+2\right)
\ - \nn\\  &&-\
\ha \left(\ervv {24} X^+_3
X^+_2 X^+_1\right) \left(a_3+2 a_2 \left(a_2+a_3+2\right)+2\right)
\ + \nn\\  &&+\ \left(a_2+1\right) \left(a_2+a_3+2\right)
\Big( 2 \left(Y_4 Y_3 X^+_{13}
X^+_2\right)  -\left(\ervv
{34} X^+_{13} X^+_2\right)
\ - \nn\\  &&-\ 2 \left(Y_4 Y_3 X^+_{23} X^+_{12}\right)
 +\left(\ervv {34} X^+_{23}
X^+_{12}\right)  +2
\left(Y_4 Y_2 X^+_3 X^+_{12}\right)  -\left(\ervv {24} X^+_3 X^+_{12}\right)
\ - \nn\\  &&-\ 2 \left(Y_4 \orv 1 X^+_3
X^+_2 \right)  +\left(\ervv
{14} X^+_3 X^+_2 \right)\Big)
+a_2 \left(a_2+1\right) \left(a_2+a_3+2\right)
\ \times \nn\\  && \times\
\Big(-4\left(Y_4 Y_2 X^+_{13}\right)  +2 \left(\ervv {24} X^+_{13}\right)  +4 \left(Y_4 \orv 1
X^+_{23}\right)  -2\left(\ervv {14} X^+_{23}\right)\Big)
\ + \nn\\  &&+\
 \left(a_2+1\right) \left(a_2+a_3+1\right) \left(a_2+a_3+2\right)
\ \times \nn\\  && \times\
\Big(-4
\left(Y_3 Y_2 X^+_{12}\right)  +2 \left(\ervv {23}
X^+_{12}\right)  +4 \left(Y_3 \orv 1 X^+_2 \right)
 -2
\left(\ervv {13} X^+_2 \right)   \ - \nn\\  &&-\  8 \left(Y_2 \orv
1\right) a_2
+4
a_2
\ervv {12}\Big) \Big)\lwv
\eea 
where the root vector ~$X^+_{jk}$~ corresponds to the compact root ~$\d_j-\d_{k+1}=\a_j+\a_{j+1}+\cdots+\a_k\,$,
 ~$Y_{k}$~ corresponds to the odd (noncompact) root ~$\d_k = \a_k+\a_{k+1}+\cdots+\a_n\,$, (thus ~$Y_4\equiv X^+_4$),
 ~$Y_{jk}$~ corresponds to the even noncompact root ~$\d_j+\d_k\,$.
In \eqref{vs12} it is more convenient to use a PBW type of basis with the compact roots ~$X^+_{...}$~ to the right
of the noncompact roots ~$Y_{...}\,$. ~The norm of \eqref{vs12} is:
\eqnn{vs12n}&& 64 a_2 \left(a_2+1\right){}^2 \left(a_2+2\right)
\left(a_2+a_3+1\right) \left(a_2+a_3+2\right){}^2
\left(a_2+a_3+3\right)\ \times \nn\\ &&\times\ \left(-2 d+a_2+a_3+4\right) \left(-2
d+a_2+a_3+5\right)/\left(2 a_2+a_3+2\right){}^2. \nn\eea 
For ~$d ~=~ d_{12}\,$, ~$a_1 = 0,\ a_2\neq 0$~ the singular vector
\eqref{vs12} is non-trivial and gives rise to a invariant subspace
which must be factored out for unitarity. For ~$d< \frh + \half
(a_2+a_3) $~ the vector \eqref{vs12} is not singular but has
negative norm and there is no unitarity for ~$a_2\neq 0$, except at
the isolated unitary point ~$d ~=~ 2 + \half (a_2+a_3) ~=~ d_{2} >
d_{13}\,$, where the vector \eqref{vs12} has zero norm and can not
spoil the unitarity. For that value of ~$d$~ there is a singular
vector ~$v^1_{\d_2}$~ of weight ~$\d_2 =\a_2 + \a_3 +\a_4$~
\cite{Dos,DZ}:
\eqnn{svn3d2} v^1_{\d_2} \, &=&\,  \sum_{k_1=0}^1
\sum_{k_2=0}^1  \, b_{k_1,k_2}\, (X^+_2)^{1-k_1} (X^+_3)^{1-k_2} \,
\times \\ &&\times\ X^+_4\, (X^+_{3})^{k_{2}}  \, (X^+_{2})^{k_{1}}
\,
 v_0\, \equiv\, \cp^{1,\d_2}\ v_0 \ ,
\nn\\ \nn\\
 b_{k_1,k_2} \, &=&\,  (-1)^{k_1 + k_2}
\,   \frac{a_2+k_1}{1+a_2+a_3 - k_2} \nn\eea
which has to be factored out for unitarity for ~$a_2\neq 0$~, while for
~$a_2=0$~ it is descendant of the compact vector ~$X_2^+\,v_0\,$.\\
Overall no further unitarity is possible for ~$a_2\neq0$, thus below we consider only the cases ~$a_1 = a_2= 0$.
Then the singular vectors above
are descendants of compact root singular vectors
~$X^+_1\,v_0\,$~ and ~$X^+_2\,v_0\,$, thus, there is no obstacle for
unitarity for ~$ d > 2 + \half a_3 ~=~ d_2 = d_{13}$ (for $a_1 = a_2= 0$).
The next reducibility point is ~$d=d_{13}=d_2$.
The singular vector for ~$d=d_{13}$~ and ~$m=1$~ has weight ~$\d_1+\d_3=\a_1+\a_2+2\a_3+2\a_4$~:
\eqnn{d1d3} v^1_{\d_1+\d_3} &=&\Big( -4 a_1 \left(Y_4 Y_3 X^+_3 X^+_{12}\right)-2
a_1 \left( Y^2_4 (X^+_3 )^2 X^+_{12}\right)- \nn\\ && -
2 \left(a_1+a_2+1\right)
\left( Y^2_4 X^+_3 X^+_{23} X^+_1\right)+4 a_1
\left(a_1+a_2+1\right) \left( Y^2_4 X^+_3 X^+_{13}\right)+ \nn\\ && +4
\left(a_3+1\right) \left(Y_4 Y_3 X^+_{23} X^+_1\right)-8 a_1
\left(a_3+1\right) \left(Y_4 Y_3 X^+_{13}\right)- \nn\\ && -
4
\left(a_1+a_2+a_3+2\right) \left(Y_4 Y_2 X^+_3 X^+_1\right)+8 a_1
\left(a_1+a_2+a_3+2\right) \left(Y_4 \orv 1 X^+_3 \right)+ \nn\\ && +8 a_3
\left(a_1+a_2+a_3+2\right) \left(Y_3 Y_2 X^+_1\right)+4 a_3
\left(a_1+a_2+a_3+2\right) \left(Y^2_3 X^+_2 X^+_1\right)-\nn\\ && -8 a_1 a_3
\left(a_1+a_2+a_3+2\right) \left(Y^2_3 X^+_{12}\right)+ \nn\\ && +2 \left(a_1
\left(a_3-1\right)+  a_2 \left(a_3-1\right)-2\right) \left(\ervv {34}
X^+_{23} X^+_1\right)-\nn\\ && -4 a_1 \left(a_1 \left(a_3-1\right)+ a_2
\left(a_3-1\right)-2\right) \left(\ervv {34} X^+_{13}\right)+\nn\\ && +2
\left(a_1+a_2+2\right) \left(a_1+a_2+a_3+2\right) \left(\ervv {24}
X^+_3 X^+_1\right)-\nn\\ && -4 a_1 \left(a_1+a_2+2\right)
\left(a_1+a_2+a_3+2\right) \left(\ervv {14} X^+_3 \right)-\nn\\ && -4
\left(a_1+a_2+2\right) a_3 \left(a_1+a_2+a_3+2\right) \left(\ervv
{23} X^+_1\right)-  \nn\\ && -\left(a_1+a_2+2 a_3+2\right) \left(\ervv {34}
X^+_3 X^+_2 X^+_1\right)+2 a_1 \left(a_1+a_2+2 a_3+2\right)
\left(\ervv {34} X^+_3 X^+_{12}\right)+\nn\\ && +8 a_1 \left(a_1+a_2+2\right)
a_3 \left(a_1+a_2+a_3+2\right) \ervv {13}-16 a_1 a_3
\left(a_1+a_2+a_3+2\right) \left(Y_3 \orv 1\right)+\nn\\ && +2 \left(Y_4 Y_3
X^+_3 X^+_2 X^+_1\right)+ Y^2_4 (X^+_3 )^2 X^+_2 X^+_1\Big)
\lwv
\nn\eea 
 For ~$a_1=a_2=0$~ it is descendant of the compact root singular vector ~$X^+_1\, v_0\,$.
  However, there is a subsingular vector:
\eqnn{vssd213} v_{2,13}^{ss} ~&=&~ \Big(
2 a_3 \left(\ervv {23} \orv 1\right)-2 a_3 \left(\ervv {13} Y_2\right)+2 a_3 \left(Y_3 \left(\ervv {12}\right)\right)-4 a_3
\left(Y_3 Y_2 \orv 1\right)\ + \nn\\ && +\ 2 \left(Y_4 Y_3 Y_2
X^+_{13}\right) - \ervv {34} Y_2 X^+_{13}+\ervv {24} Y_3
X^+_{13}-Y_4 \ervv {23} X^+_{13}-2 \left(Y_4 Y_3 \orv 1
X^+_{23}\right) \ + \nn\\ && +\
\ervv {34} \orv 1 X^+_{23}-\ervv {14} Y_3
X^+_{23}+Y_4 \ervv {13} X^+_{23}+2 \left(Y_4 Y_2 \orv 1
X^+_3 \right)-\ervv {24} \orv 1 X^+_3
\ + \nn\\ && +\
\ervv {14} Y_2 X^+_3 -Y_4 \left(\ervv {12}\right) X^+_3 \Big) \lwv
 \eea 
with norm:
$$ -16 a_3 \left(a_3+3\right) \left(-2 d+a_3+2\right)
\left(-2 d+a_3+3\right) \left(-2 d+a_3+4\right) $$
This vector must
be factorized in order to obtain UR at ~$d=d_2 = d_{13}\,$. But
below this value the vector \eqref{vssd213} has negative norm if
~$a_3\neq 0$~ and there is no unitarity, except at the isolated
unitary point ~$d ~=~ \trh + \half a_3 ~=~ d_{23} > d_{14}\,$. At
that value of $d$ there is a singular vector of weight ~$\d_2+\d_3 =
\a_2 + 2\a_3$~:
 \eqnn{svd23} v^1_{\d_2+\d_3} &=& \Big(
2 \left(a_3+1\right) \left(Y_4 Y_3 X^+_{23}\right)-2
\left(a_3+1\right) \left(Y_4 Y_2 X^+_3 \right)+2 a_3
\left(a_3+1\right) \left(Y^2_3 X^+_2
\right)-\nn\\ && -\left(a_3+1\right) \left(\ervv {34}
X^+_{23}\right)+\left(a_3+1\right) \left(\ervv {24} X^+_3
\right)-\frac{1}{2} \left(2 a_3+1\right) \left(\ervv {34} X^+_3
X^+_2 \right)-\\ && -2 a_3 \left(a_3+1\right) \ervv {23}+4 a_3
\left(a_3+1\right) \left(Y_3 Y_2\right)+Y_4 Y_3 X^+_3
X^+_2 +\frac{1}{2} \left( Y^2_4 (X^+_3 )^2 X^+_2 \right)
\Big) \lwv
\nn\eea 
with norm:
$$ 16 a_3 \left(a_3+1\right){}^2 \left(a_3+2\right) \left(-2 d+a_3+2\right) \left(-2 d+a_3+3\right)$$
For ~$a_3\neq 0$~ the singular vector \eqref{svd23} should be
factored for unitarity, while for ~$a_3=0$~ it is descendant of the
compact singular vectors.\\
In the same range for ~$a_3\neq 0$~ at ~$d=d_3 =1+\half a_3$~
there is a singular vector of weight ~$\d_3 = \a_3+\a_4$~:
\eqn{svd3}
v^1_{\d_3}\, =\,  \sum_{k=0}^1   \,(-1)^{k}\, (a_3+k) \,
(X^+_3)^{1-k}    \, X^+_4\, (X^+_{3})^{k}  \,
 v_0\, \equiv\, \cp^{1,\d_3}\, v_0
 \ee
 which must be  factored out for unitarity.\\
On the other hand, for ~$a_1 = a_2= a_3=0$~ all (sub)singular vectors above
are descendants of the compact singular vectors ~$X^+_k\,v_0\,$,
$k=1,2,3$,
 and there is no obstacle for unitarity for ~$ d > \trh ~=~ d_{23} = d_{14}$.
For ~$a_3=0$~ and ~$d ~=~ \trh$~ there is also a singular vector of weight ~$\d_1+\d_4$~:
\eqnn{} v^1_{\d_1+\d_4} ~&=&~ \Big(-4 \left(Y_4 Y_2 X^+_1\right)-
2 \left( Y^2_4 X^+_{23} X^+_1\right)+2 \left(Y_4 Y_3 X^+_2 X^+_1\right)+ \nn\\ && +
Y^2_4 X^+_3 X^+_2 X^+_1- 3 \left(\ervv {34} X^+_2 X^+_1\right)+6 \left(\ervv {24} X^+_1\right)\Big)\lwv
\nn\eea 
but it is also descendant of compact singular vectors. Finally, for
~$d ~=~ \trh$~ there is a subsingular vector of weight
~$\d_1+\d_2+\d_3+\d_4$~:
\eqn{ss2314} v^{ss}_{\d_1+\d_2+\d_3+\d_4}
~=~ \sum_{i,j,k,\ell=1}^4 ~\eps^{ijk\ell}\ Y_i\, Y_j\,
Y_k\, Y_\ell\, \ v_0 \ee where ~$\eps^{ijk\ell}$~ is the
totally antisymmetric symbol so that ~$\eps^{1234}=1$. The norm of
the vector \eqref{ss2314} is ~~$2304(-1 + d)d(-3 + 2d)(-1 + 2d)$.
Thus, for ~$\trh > d >1$~ there is no unitarity since then the vector \eqref{ss2314}
has negative norm. In all cases there will be no unitarity for ~$d\leq 1$~ except possibly when
~$a_1 = a_2= a_3=0$~ to which we restrict below.
At ~$d=d_3 =d_{24}=1$~ there are the singular vector \eqref{svd3} and
the singular vector  of weight ~$\d_2+\d_4 = \a_2 + \a_3 + 2\a_4\,$~:
\eqnn{svd24} v^1_{\d_2+\d_4}
~=~ \Big( -\left(a_3+2\right) \left(\ervv {34} X^+_2 \right)+2
\left(Y_4 Y_3 X^+_2 \right)+ Y^2_4 X^+_3 X^+_2 \Big)\lwv
\eea 
both of which are descendants of compact singular vectors. At ~$d=d_3 =d_{24}=1$~
there is also a subsingular vector \be v^{ss}_{\d_2+\d_3+\d_4} ~=~
\left( Y_2\, Y_3\, Y_4 - Y_4\, Y_3\, Y_2 \right) ~=~ \third
\sum_{i,j,k=2}^4 ~\eps^{ijk}\ Y_i\, Y_j\, Y_k\, \ v_0 \ees of norm
~~$144 d(d-1)(2d-1)$. It is not an obstacle for unitarity for $d=1$,
but for ~$d<1$. Thus, there is no unitarity for ~$d<1$ except  at
the isolated unitary point ~$d = \half = d_{34}$. At that point all
(sub)singular vectors above are descendants of compact singular
vectors. Yet there is  the singular vector:
\be v^1_{\d_3+\d_4} ~=~
\left( \half Y^2_4\,X^+_3 - 2 Y_3\,Y_4 + Y_{34} \right) \, v_0
 \ees
with norm ~~$8d(2d-1)$. It is not an obstacle for unitarity for $d=\half$, but for ~$d<\half$.
Thus, there is no unitarity for ~$d<\half$ except at the isolated point
 ~$d=d_4=0=a_1 = a_2=a_3$~
where we have the trivial one-dimensional UIR since all possible states are descendants of
factored out singular vectors.\dia

\end{document}